\documentclass[12pt]{amsart}

\usepackage{aliascnt,amsmath,amssymb,amsthm,amsfonts,enumerate,mathrsfs,latexsym,mathtools,bm,xcolor,tikz}
\usepackage[abbrev]{amsrefs}

\usepackage[OT2,T1]{fontenc}
\usepackage{scalerel}

\usepackage[marginparwidth=0pt,margin=24truemm]{geometry}

\definecolor{mylinkcolor}{RGB}{16, 156, 81}
\definecolor{mycitecolor}{RGB}{20, 80, 140}

\usepackage{hyperref}
\usepackage[nameinlink]{cleveref}
\hypersetup{
    setpagesize=false, 
    bookmarksnumbered=true, 
    bookmarksopen=true, 
    colorlinks=true, 
    linkcolor=mylinkcolor, 
    citecolor=mycitecolor,
}
\numberwithin{equation}{section}

\theoremstyle{plain}

\newtheorem{thm}{Theorem}[section]
\crefname{thm}{Theorem}{Theorems}

\newaliascnt{lem}{thm} 
\newtheorem{lem}[lem]{Lemma}
\aliascntresetthe{lem}
\crefname{lem}{Lemma}{Lemmas}

\newaliascnt{prop}{thm}
\newtheorem{prop}[prop]{Proposition}
\aliascntresetthe{prop}
\crefname{prop}{Proposition}{Propositions}

\newtheorem*{mainthm}{Main Theorem}

\allowdisplaybreaks[2]

\everymath{\displaystyle}

\newcommand{\summ}[1]{\sum_{\substack{#1}}}

\newcommand{\GG}{G}
\newcommand{\GD}{\operatorname{D}}

\newcommand{\QQ}{\mathbb{Q}}
\newcommand{\ZZ}{\mathbb{Z}}
\newcommand{\RR}{\mathbb{R}}
\newcommand{\CC}{\mathbb{C}}
\newcommand{\HH}{\mathbb{H}}

\DeclareMathOperator{\id}{id}

\renewcommand{\Im}{\operatorname{Im}}

\title{Zagier's $2$-$3$-$2$ formula for multiple Eisenstein series}

\author{Henrik Bachmann}
\address{Graduate School of Mathematics, Nagoya University, Nagoya, Japan.}
\email{henrik.bachmann@math.nagoya-u.ac.jp}

\date{\today}


\begin{document}

\begin{abstract}
We prove a multiple Eisenstein series analogue of Zagier's formula for the multiple zeta values $\zeta(2,\ldots,2,3,2,\ldots,2)$. The multiple Eisenstein series identity contains correction terms given by derivatives of odd Eisenstein series. Using the Goncharov coproduct and Brown's motivic version of Zagier's formula, we show that the normalized difference has rational Fourier coefficients. We then use the $\mathfrak{sl}_2$-structure and a modularity criterion of Yu to prove the identity by induction.
\end{abstract}

\maketitle

\section{Introduction}

The purpose of this paper is to prove a multiple Eisenstein series analogue of Zagier's formula for the multiple zeta values $\zeta(2,\ldots,2,3,2,\ldots,2)$ \cite{Z}, where
\[
 \zeta(k_1,\ldots,k_r)=\sum_{m_1>\cdots>m_r>0}\frac{1}{m_1^{k_1}\cdots m_r^{k_r}}
 \qquad(k_1\geq2,\ k_2,\ldots,k_r\geq1).
\]
We write $\{2\}^m$ for a string of $m$ copies of $2$, with $\{2\}^0$ the empty string and $\zeta(\{2\}^0)=1$. Zagier's evaluation expresses $\zeta(\{2\}^a,3,\{2\}^b)$ as a rational linear combination of products of an odd zeta value and $\zeta(2,\ldots,2)$, see \eqref{eq:Zagier-MZV} below. This formula was important for Brown's proof of Hoffman's conjecture that every multiple zeta value is a $\QQ$-linear combination of those whose entries are two or three \cite{Br}. Let $\tau\in\HH=\{z\in\CC\mid\Im(z)>0\}$ and $\ZZ_M=\{m\in\ZZ\mid |m|<M\}$. For $k_1,\ldots,k_r\geq2$, the multiple Eisenstein series introduced by Gangl, Kaneko and Zagier \cite{GKZ} are
\[
 \GG(k_1,\ldots,k_r;\tau)=\lim_{M\to\infty}\lim_{N\to\infty}\summ{\lambda_1,\ldots,\lambda_r\in\ZZ_M\tau+\ZZ_N\\ \lambda_1\succ\cdots\succ\lambda_r\succ0}\frac{1}{\lambda_1^{k_1}\cdots\lambda_r^{k_r}}\,.
\]
Here $\lambda_1\succ\lambda_2$ means $\lambda_1-\lambda_2\in\HH\cup\RR_{>0}$. This limit is absolutely convergent if $k_1\geq3$. Multiple Eisenstein series can be viewed as hybrids of multiple zeta values and modular forms. Write $\mathcal Z$ for the $\QQ$-algebra of multiple zeta values and put $q=e^{2\pi i\tau}$. The Fourier expansions of multiple Eisenstein series have the form
\[
 \GG(k_1,\ldots,k_r;\tau)=\zeta(k_1,\ldots,k_r)+\sum_{n\geq1}a_nq^n,
 \qquad a_n\in\mathcal Z[2\pi i]\,.
\]
A general question is which relations among multiple zeta values lift to multiple Eisenstein series. Some relations lift without correction, but not all. Derivatives are not visible at the constant term, and in \cite{BKM} the author, Kanno and Maesaka proved that the $\QQ$-algebra of multiple Eisenstein series is closed under derivatives. Relations involving cusp forms give another source, since cusp forms have vanishing constant term. The same work gives a family of relations and conjectures that these relations together with their derivatives generate all linear relations among multiple Eisenstein series. Here we show that Zagier's relation lifts after adding correction terms given by derivatives of odd Eisenstein series multiplied by the quasimodular forms $\GG(2,\ldots,2)$ studied in \cite{Ba}. In the following, we usually leave out the variable and write $\GG(k_1,\ldots,k_r)=\GG(k_1,\ldots,k_r;\tau)$, with $\GG(\{2\}^0)=1$. Put
\[
 \GD=2\pi i\frac{d}{d\tau}=(2\pi i)^2q\frac{d}{dq}\,.
\]
This derivative raises the weight by two. Let $\mathcal E$ be the $\QQ$-linear space spanned by $1$ and the multiple Eisenstein series $\GG(k_1,\ldots,k_r)$ with $r\geq1$ and $k_1,\ldots,k_r\geq2$. The harmonic product makes $\mathcal E$ a $\QQ$-algebra, and \cite{BKM}*{Main Theorem~D} shows that it is closed under $\GD$. In \cite{BK}, the author and Kanno show that $\mathcal E$ is graded by weight and has a derivation $\delta$ which lowers the weight by two and satisfies
\[
 \delta\GG(k_1,\ldots,k_r)=
 \begin{cases}
  -\frac12\GG(k_2,\ldots,k_r),&k_1=2,\\
  0,&k_1\geq3,
 \end{cases}
\]
where $\GG(\varnothing)=1$. For $F\in\mathcal E$ of weight $k$, we have
\begin{equation}\label{eq:lowering-commutator}
 \delta\GD F-\GD\delta F=kF\,.
\end{equation}
Together with the weight operator, $\GD$ and $\delta$ give $\mathcal E$ an $\mathfrak{sl}_2$-algebra structure. We use the lowering operator to prove the formula below by induction on the number of initial entries equal to two.

For $a,b\geq0$ and $j\geq1$, put
\[
 c_{a,b}^{j,0}=2(-1)^j\left(
 \binom{2j}{2b+2}-\bigl(1-4^{-j}\bigr)\binom{2j}{2a+1}
 \right).
\]
These are the coefficients in Zagier's formula \cite{Z},
\begin{equation}\label{eq:Zagier-MZV}
 \zeta(\{2\}^a,3,\{2\}^b)=\sum_{j=1}^{a+b+1}c_{a,b}^{j,0}\zeta(\{2\}^{a+b+1-j})\zeta(2j+1)\,.
\end{equation}
Set $c_{a,b}^{j,m}=0$ when $a<0$, and for $m\geq0$ define recursively
\begin{equation}\label{eq:zagier-coefficient-recurrence}
 c_{a,b}^{j,m+1}=\frac{c_{a,b}^{j,m}-c_{a-1,b}^{j,m}}{2(m+1)(2j+m+1)}.
\end{equation}
\begin{mainthm}
For $a,b\geq0$, we have
\begin{equation}\label{eq:zagier-MES}
 \GG(\{2\}^a,3,\{2\}^b)
 =\sum_{\substack{j\geq1,\ m\geq0\\j+m\leq a+b+1}}
 c_{a,b}^{j,m}\GG(\{2\}^{a+b+1-j-m})
 \GD^m\!\GG(2j+1).
\end{equation}
\end{mainthm}
\noindent We explain the proof for $a=b=1$, where \eqref{eq:zagier-MES} gives
\begin{equation}\label{eq:232}
 \GG(2,3,2)=-\frac{11}{2}\GG(2)\GG(5)+\frac{75}{8}\GG(7)-\frac12\GG(2)\GD\GG(3)-\frac1{16}\GD^2\!\GG(3)\,.
\end{equation}
The lowering operator relates this identity to the corresponding formula for $\GG(3,2)$. In \cref{prop:normalized-reality}, we use the Goncharov coproduct and Brown's motivic version of Zagier's formula to show that the difference between the two sides of \eqref{eq:232} has rational Fourier coefficients after division by $(-2\pi i)^7$. The same result in weight five gives
\[
 (-2\pi i)^{-5}\left(\GG(3,2)-3\GG(2)\GG(3)+\frac{11}{2}\GG(5)-\frac12\GD\GG(3)\right)\in\QQ[[q]]\,.
\]
The commutator relation gives $\delta\GD\GG(3)=3\GG(3)$ and $\delta\GD^2\GG(3)=8\GD\GG(3)$. Since $\delta$ is a derivation, we have
\[
 \delta\left(\GG(3,2)-3\GG(2)\GG(3)+\frac{11}{2}\GG(5)-\frac12\GD\GG(3)\right)=0,
\]
and
\begin{align*}
 &\delta\left(\GG(2,3,2)+\frac{11}{2}\GG(2)\GG(5)-\frac{75}{8}\GG(7)+\frac12\GG(2)\GD\GG(3)+\frac1{16}\GD^2\GG(3)\right)\\
 &\qquad=-\frac12\left(\GG(3,2)-3\GG(2)\GG(3)+\frac{11}{2}\GG(5)-\frac12\GD\GG(3)\right).
\end{align*}
In positive weight, \cite{BK}*{Main Theorem~A(iii)} identifies the kernel of $\delta$ with the span of the multiple Eisenstein series whose first entry is at least three. The weight-five difference above lies in this kernel and has real normalized Fourier coefficients, so Yu's theorem \cite{Y} forces it to vanish. It follows that the weight-seven difference also lies in this kernel, and Yu's theorem gives \eqref{eq:232}. The same argument gives the induction on $a$ in \Cref{sec:proof}.

\section*{Acknowledgments}
This project was partially supported by JSPS KAKENHI Grant Number JP26K22254.

\section{The odd motivic terms}\label{sec:reduction}

In this section, we show that the difference in the main theorem has rational Fourier coefficients after normalization. For $k_1,\ldots,k_r\geq2$, put
\[
 \tilde g(k_1,\ldots,k_r)=\frac{(-2\pi i)^{k_1+\cdots+k_r}}{(k_1-1)!\cdots(k_r-1)!}\summ{m_1>\cdots>m_r>0\\ n_1,\ldots,n_r>0}n_1^{k_1-1}\cdots n_r^{k_r-1}q^{m_1n_1+\cdots+m_rn_r}\,.
\]
In \cite{BT}, the author and Tasaka described the Fourier expansions of multiple Eisenstein series in terms of the Goncharov coproduct, which we combine with Brown's motivic version of Zagier's formula \cite{Br}. Using the notation of \cite{BT}, we view the family of multiple Eisenstein series as a map $\GG$ on the formal iterated-integral algebra. By \cite{BT}*{Definition~4.8 and Theorem~4.10}, this map is the algebra homomorphism given by
\begin{equation}\label{eq:Goncharov-realization}
 \GG=\mathrm m\circ\bigl(\tilde g\otimes\zeta\bigr)\circ\Delta_{\mathrm G}\,.
\end{equation}
Here $\mathrm m$ is multiplication, $\Delta_{\mathrm G}$ is the Goncharov coproduct, and $\tilde g$ and $\zeta$ denote the corresponding character extensions on the formal iterated-integral algebra of \cite{BT}.\footnote{Notice that we reverse the indices of \cite{BT}, hence the order of the blocks $10^{k_i-1}$ in the iterated-integral words. We also order the two tensor factors in \eqref{eq:Goncharov-realization} as in \cite{G}, not as in \cite{BT}.}

Let $\mathcal H^{\mathfrak m}$ denote the weight-graded $\QQ$-algebra of motivic multiple zeta values of \cite{Br}, with period map $\operatorname{per}\colon\mathcal H^{\mathfrak m}\to\CC$. Brown's motivic iterated integrals $I^{\mathfrak m}(a_0;a_1,\ldots,a_M;a_{M+1})$ with $a_0,\ldots,a_{M+1}\in\{0,1\}$ satisfy the defining relations of the formal iterated-integral algebra \cite{Br}*{Section~2}. We write $\zeta^{\mathfrak m}$ for the resulting character, which satisfies $\operatorname{per}\circ\zeta^{\mathfrak m}=\zeta$. Write $\mathcal H_q$ for the $\QQ$-algebra generated by $\GD^m\tilde g(\boldsymbol k)$ for $m\geq0$ and indices $\boldsymbol k$ with entries at least two. This algebra is closed under $\GD$ by definition. Replacing $\zeta$ in \eqref{eq:Goncharov-realization} by $\zeta^{\mathfrak m}$ gives
\begin{equation}\label{eq:motivic-lift}
 \widetilde\GG=\bigl(\tilde g\otimes\zeta^{\mathfrak m}\bigr)\circ\Delta_{\mathrm G}\,,
\end{equation}
which satisfies $\mathrm m\circ(\id\otimes\operatorname{per})\circ\widetilde\GG=\GG$. On indices with entries at least two, its values lie in $\mathcal H_q\otimes_{\QQ}\mathcal H^{\mathfrak m}$. For $a,b\geq0$, put
\begin{equation}\label{eq:Xi-difference}
\begin{aligned}
 \Xi_{a,b}={}&\widetilde\GG(\{2\}^a,3,\{2\}^b)
 -\sum_{\substack{j\geq1,\ m\geq0\\ j+m\leq a+b+1}}c_{a,b}^{j,m}\,\widetilde\GG(\{2\}^{a+b+1-j-m})\\
 &\quad\cdot\left(\GD^m\tilde g(2j+1)\otimes1+\boldsymbol1_{m=0}\,1\otimes\zeta^{\mathfrak m}(2j+1)\right)
 \in\mathcal H_q\otimes_{\QQ}\mathcal H^{\mathfrak m}.
\end{aligned}
\end{equation}
Since $\GD$ kills constants and $\GG(2j+1)=\zeta(2j+1)+\tilde g(2j+1)$, applying $\mathrm m\circ(\id\otimes\operatorname{per})$ to $\Xi_{a,b}$ gives the difference of the two sides of \eqref{eq:zagier-MES}.

The infinitesimal coaction operators of \cite{Br} act on the second factor. Let $\mathcal L=\textstyle\bigoplus_N\mathcal L_N$ be Brown's graded Lie coalgebra of indecomposables modulo $\zeta^{\mathfrak m}(2)$, and write $\operatorname{pr}\colon\mathcal H^{\mathfrak m}_{>0}\to\mathcal L$ for the quotient map. For odd $N\geq3$, the derivation $D_N\colon\mathcal H^{\mathfrak m}\to\mathcal L_N\otimes\mathcal H^{\mathfrak m}$ of \cite{Br}*{Definition 3.1} is given by
\begin{align*}
 D_NI^{\mathfrak m}(a_0;a_1,\ldots,a_M;a_{M+1})={}&\sum_{p=0}^{M-N}\operatorname{pr}\bigl(I^{\mathfrak m}(a_p;a_{p+1},\ldots,a_{p+N};a_{p+N+1})\bigr)\\
 &\hspace{18mm}\otimes I^{\mathfrak m}(a_0;a_1,\ldots,a_p,a_{p+N+1},\ldots,a_M;a_{M+1})\,.
\end{align*}
Brown's kernel theorem \cite{Br}*{Theorem 3.3}, Euler's relation and \cite{Br}*{Lemma 3.2} give
\[
 \bigcap_{r\geq1}\ker D_{2r+1}=\QQ[\zeta^{\mathfrak m}(2)],
 \qquad
 D_{2r+1}\zeta^{\mathfrak m}(2j+1)
 =\boldsymbol1_{r=j}\,\operatorname{pr}\bigl(\zeta^{\mathfrak m}(2r+1)\bigr)\otimes1.
\]
On $\mathcal H_q\otimes\mathcal H^{\mathfrak m}$, we write $\mathcal D_N=\id\otimes D_N$, with the $\mathcal L_N$-factor moved to the front.

For the cut calculation, we associate to an index $(k_1,\ldots,k_r)$ with all entries at least two and weight $K$ the word $\hat w(k_1,\ldots,k_r)=10^{k_r-1}10^{k_{r-1}-1}\cdots10^{k_1-1}$, placed between the base points of $I(0;\hat w;1)$, so that $\zeta^{\mathfrak m}(k_1,\ldots,k_r)=I^{\mathfrak m}(0;\hat w(k_1,\ldots,k_r);1)$. In particular
\[
 \hat w(\{2\}^u)=(10)^u,
 \qquad
 \hat w(\{2\}^a,3,\{2\}^b)=(10)^b100(10)^a\,.
\]
We write $\hat a_1\cdots\hat a_K$ for the letters of $\hat w$ and extend by the base points $\hat a_0=0$ and $\hat a_{K+1}=1$. The extended word of $(\{2\}^a,3,\{2\}^b)$ is alternating except for the central pair $\hat a_{2b+2}=\hat a_{2b+3}=0$. For these words the Goncharov coproduct \cite{G} is
\begin{equation}\label{eq:Goncharov-cuts}
 \Delta_{\mathrm G}\,\hat w=\sum_{0=i_0<i_1<\cdots<i_k<i_{k+1}=K+1}(\hat a_{i_1}\cdots\hat a_{i_k})\otimes\prod_{p=0}^{k}I(\hat a_{i_p};\hat a_{i_p+1},\ldots,\hat a_{i_{p+1}-1};\hat a_{i_{p+1}})\,.
\end{equation}
Each choice $i_1<\cdots<i_k$ in \eqref{eq:Goncharov-cuts} is a cut. The selected letters $\hat a_{i_1}\cdots\hat a_{i_k}$ form its kept word in the first tensor factor, while the iterated integrals between consecutive selected points are the pieces in the second factor.

\begin{figure}[htbp]
\centering
\begin{tikzpicture}[scale=.85]
 \draw[gray] (180:2) arc[start angle=180,end angle=0,radius=2];
 \draw[thick] (180:2)--(112.5:2)--(90:2)--(67.5:2)--(0:2)--cycle;
 \foreach \t in {180,157.5,135,112.5,90,67.5,45,22.5,0}
  \filldraw[fill=white] (\t:2) circle (1.6pt);
 \foreach \t in {180,112.5,90,67.5,0}
  \fill (\t:2) circle (2.2pt);
 \node at (180:2.25) {$0$};
 \node at (157.5:2.25) {$1$};
 \node at (135:2.25) {$0$};
 \node at (112.5:2.25) {$1$};
 \node at (90:2.25) {$0$};
 \node at (67.5:2.25) {$0$};
 \node at (45:2.25) {$1$};
 \node at (22.5:2.25) {$0$};
 \node at (0:2.25) {$1$};
\end{tikzpicture}

\normalfont\small The cut of $\hat w(2,3,2)=1010010$ shown here contributes one summand to \eqref{eq:Goncharov-cuts}. Its filled inner vertices give the kept word $100$, and the arcs between consecutive filled vertices represent the pieces.
\end{figure}

A nonempty piece with equal boundary points vanishes (see \cite{G}). The boundary pieces therefore force a nonempty kept word to start with $1$ and end with $0$. Since the input word has no consecutive $1$'s, consecutive kept $1$'s would bound a nonempty piece with equal endpoints. Hence a nonempty kept word is $10^{l_r-1}\cdots10^{l_1-1}$ with $l_1,\ldots,l_r\geq2$ and evaluates as $\tilde g(l_1,\ldots,l_r)$ in \eqref{eq:motivic-lift}. The empty kept word contributes $1$.

\begin{lem}
For every index $\boldsymbol k$ with entries at least two and word $\hat w(\boldsymbol k)=\hat a_1\cdots\hat a_K$, and every odd $N\geq3$, one has
\begin{equation}\label{eq:lifted-coaction}
 \mathcal D_N\,\widetilde\GG(\boldsymbol k)=\sum_{p=0}^{K-N}\operatorname{pr}\bigl(I^{\mathfrak m}(\hat a_p;\hat a_{p+1},\ldots,\hat a_{p+N};\hat a_{p+N+1})\bigr)\otimes\widetilde\GG\bigl(\hat a_1\cdots\hat a_p\,\hat a_{p+N+1}\cdots\hat a_K\bigr)\,,
\end{equation}
where $\widetilde\GG$ of a word is defined by applying \eqref{eq:motivic-lift} to \eqref{eq:Goncharov-cuts}.
\end{lem}

\begin{proof}
Apply $\mathcal D_N$ to \eqref{eq:motivic-lift}. Since $D_N$ is a derivation, each term consists of a cut and a block of weight $N$ contained in one of its pieces. For a fixed block $(\hat a_p;\hat a_{p+1},\ldots,\hat a_{p+N};\hat a_{p+N+1})$, deleting its interior identifies these cuts with the cuts of the shortened word $\hat a_1\cdots\hat a_p\hat a_{p+N+1}\cdots\hat a_K$. Summing over $p$ gives \eqref{eq:lifted-coaction}.
\end{proof}

\begin{thm}\label{thm:odd-cancellation}
For all $a,b\geq0$ and all $r\geq1$, one has $\mathcal D_{2r+1}\Xi_{a,b}=0$. In particular,
\[
 \Xi_{a,b}\in\mathcal H_q\otimes_{\QQ}\QQ[\zeta^{\mathfrak m}(2)]\,.
\]
\end{thm}

\begin{proof}
Put $n=a+b$. Since the extended word $0(10)^u1$ is alternating, every odd block in \eqref{eq:lifted-coaction} has equal boundary letters. Hence
\[
 \mathcal D_N\widetilde\GG(\{2\}^u)=0.
\]

Let $1\leq r\leq n+1$ and $N=2r+1$. For $(\{2\}^a,3,\{2\}^b)$, only the blocks crossing the central pair $00$ survive, and deleting their interiors gives $0(10)^{n+1-r}1$. Thus
\[
 \mathcal D_{2r+1}\widetilde\GG(\{2\}^a,3,\{2\}^b)
 =\operatorname{pr}\left(
 \sum_{p=0}^{2n+3-N}
 I^{\mathfrak m}(\hat a_p;\hat a_{p+1},\ldots,\hat a_{p+N};\hat a_{p+N+1})
 \right)\otimes\widetilde\GG(\{2\}^{n+1-r}).
\]
The sum in its first factor is the first factor of $D_{2r+1}\zeta^{\mathfrak m}(\{2\}^a,3,\{2\}^b)$. Brown's motivic version of Zagier's formula \cite{Br}*{Theorem 4.3}, applied with its two parameters interchanged, therefore identifies it with $c^{r,0}_{a,b}\operatorname{pr}(\zeta^{\mathfrak m}(2r+1))$.

After applying $\mathcal D_{2r+1}$ to \eqref{eq:Xi-difference}, the factors $\widetilde\GG(\{2\}^u)$ give zero, while
\[
 \mathcal D_{2r+1}\left(\GD^m\tilde g(2j+1)\otimes1+\boldsymbol1_{m=0}\,1\otimes\zeta^{\mathfrak m}(2j+1)\right)
 =\boldsymbol1_{m=0}\boldsymbol1_{r=j}\,
 \operatorname{pr}\bigl(\zeta^{\mathfrak m}(2r+1)\bigr)\otimes1.
\]
The coaction of the first term in \eqref{eq:Xi-difference} therefore gives exactly Zagier's coefficient $c^{r,0}_{a,b}$. For $m>0$, the derivative factor has motivic second factor $1$, while all odd coactions of the accompanying all-two factor vanish. The terms with $m=0$ cancel, and the case $r>n+1$ follows from weights.

Finally, write $\Xi_{a,b}=\textstyle\sum_ie_i\otimes h_i$ with the $e_i$ linearly independent. Every $h_i$ is killed by all $D_{2r+1}$, so the kernel statement above gives $h_i\in\QQ[\zeta^{\mathfrak m}(2)]$.
\end{proof}

\begin{samepage}
For $a,b\geq0$, put
\begin{equation}\label{eq:U-difference}
 U_{a,b}=\GG(\{2\}^a,3,\{2\}^b)-\sum_{\substack{j\geq1,\ m\geq0\\j+m\leq a+b+1}}c_{a,b}^{j,m}\GG(\{2\}^{a+b+1-j-m})\GD^m\GG(2j+1)\,.
\end{equation}

\begin{prop}\label{prop:normalized-reality}
For all $a,b\geq0$, we have
\[
 (-2\pi i)^{-2a-2b-3}U_{a,b}\in\QQ[[q]]\,.
\]
\end{prop}
\end{samepage}

\begin{proof}
Put $w=2a+2b+3$. In motivic weight $d$, every coefficient of $\Xi_{a,b}$ belongs to $(-2\pi i)^{w-d}\QQ[[q]]$. To see this, a cut whose pieces have total weight $d$ keeps $w-d$ letters, and $\GD=(-2\pi i)^2q\,d/dq$. The same weight count applies to each product in \eqref{eq:Xi-difference}. By \cref{thm:odd-cancellation}, only powers of $\zeta^{\mathfrak m}(2)$ remain. Choosing a homogeneous rational basis of the motivic algebra containing these powers, we can therefore write
\[
 \Xi_{a,b}=\sum_{s=0}^{a+b+1}(-2\pi i)^{w-2s}f_s(q)\otimes\zeta^{\mathfrak m}(2)^s,
 \qquad f_s(q)\in\QQ[[q]]\,.
\]
Taking periods and using $\zeta(2)=-(-2\pi i)^2/24$ gives
\[
 (-2\pi i)^{-w}U_{a,b}=\sum_{s=0}^{a+b+1}\left(-\frac1{24}\right)^sf_s(q)\in\QQ[[q]]\,.
\]
\end{proof}

\section{Proof of the Main Theorem}\label{sec:proof}

Let $\mathcal E^{\mathrm{adm}}$ be the $\QQ$-linear span of $1$ and the series $\GG(k_1,\ldots,k_r)$ with $k_1\geq3$ and $k_2,\ldots,k_r\geq2$. Write $\mathcal E_k^{\mathrm{adm}}$ for its weight-$k$ part. By the result of the author and Kanno \cite{BK}*{Main Theorem~A(iii)}, the derivation $\delta$ from the introduction satisfies $\ker\delta=\mathcal E^{\mathrm{adm}}$. We use the following consequence of Yu's modularity criterion \cite{Y}*{Theorem~1.3(ii)}.

\begin{prop}\label{prop:admissible-reality}
For odd $k\geq3$ and $F\in\mathcal E_k^{\mathrm{adm}}$, we have
\[
 (-2\pi i)^{-k}F\in\RR[[q]]\quad\Longrightarrow\quad F=0\,.
\]
\end{prop}

\begin{proof}
Since the order of the indices in \cite{Y} is opposite to ours, $F$ is a linear combination of admissible multiple Eisenstein series of weight $k$ in the sense of \cite{Y}. Its coefficients are rational, hence real. Since $k$ is odd, the hypothesis says that every Fourier coefficient of $F$ is purely imaginary. The cited theorem of Yu gives $F=0$.
\end{proof}

\begin{lem}\label{lem:lowering-error}
For $a,b\geq0$, we have
\[
 \delta U_{a,b}=\begin{cases}0,&a=0,\\-\frac12U_{a-1,b},&a\geq1.\end{cases}
\]
\end{lem}

\begin{proof}
Since $\delta\GG(2j+1)=0$, induction using \eqref{eq:lowering-commutator} gives
\[
 \delta\GD^m\GG(2j+1)=m(2j+m)\GD^{m-1}\GG(2j+1)
 \qquad(j,m\geq1)\,.
\]
We also have $\delta1=0$ and $\delta\GG(\{2\}^u)=-\GG(\{2\}^{u-1})/2$ for $u\geq1$. Applying $\delta$ to the right side of \eqref{eq:zagier-MES} therefore gives
\[
 \sum_{\substack{j\geq1,\ m\geq0\\j+m\leq a+b}}
 \left(-\frac12c_{a,b}^{j,m}+(m+1)(2j+m+1)c_{a,b}^{j,m+1}\right)
 \GG(\{2\}^{a+b-j-m})\GD^m\GG(2j+1)\,.
\]
By \eqref{eq:zagier-coefficient-recurrence}, the coefficient in parentheses is $-c_{a-1,b}^{j,m}/2$. For $a\geq1$, this is $-1/2$ times the right side of \eqref{eq:zagier-MES} with $a$ replaced by $a-1$. When $a=0$, every coefficient vanishes. Applying $\delta$ to the left side of \eqref{eq:zagier-MES} removes its first entry with coefficient $-1/2$ if $a\geq1$, and gives zero if $a=0$. Subtracting proves the lemma.
\end{proof}

\begin{proof}[Proof of the Main Theorem]
Fix $b\geq0$. We prove $U_{a,b}=0$ by induction on $a$. For $a=0$, \cref{lem:lowering-error} gives $\delta U_{0,b}=0$, and for $a\geq1$ the same conclusion follows from the induction hypothesis. Thus $U_{a,b}\in\mathcal E_{2a+2b+3}^{\mathrm{adm}}$. Its normalized Fourier coefficients are rational by \cref{prop:normalized-reality}, so \cref{prop:admissible-reality} gives $U_{a,b}=0$.
\end{proof}

\vspace{1cm}

{\bf AI \& computational resource disclosure:} The formula in the Main Theorem and the idea of the proof using Yu's result and the $\mathfrak{sl}_2$-action are due to the author. ChatGPT Astra and Claude Opus 5.5 worked out the proofs and checked them for correctness. The author has checked the proofs himself and takes full responsibility for the results. Computations were carried out and checked in SageMath.


\begin{thebibliography}{99}
\bibitem[Ba]{Ba} H.~Bachmann, \textit{MacMahon's sums-of-divisors and their connection to multiple Eisenstein series}, Res. Number Theory \textbf{10} (2024), Art.~50.

\bibitem[BK]{BK} H.~Bachmann and H.~Kanno, \textit{The $\mathfrak{sl}_2$-algebra structure of multiple Eisenstein series}, preprint, \href{https://arxiv.org/abs/2609.03777}{arXiv:2609.03777}, 2026.

\bibitem[BKM]{BKM} H.~Bachmann, H.~Kanno and T.~Maesaka, \textit{Relations and derivatives of multiple Eisenstein series}, preprint, \href{https://arxiv.org/abs/2602.08176v2}{arXiv:2602.08176v2}, 2026.

\bibitem[Br]{Br} F.~Brown, \textit{Mixed Tate motives over $\mathbb{Z}$}, Ann. of Math. (2) \textbf{175} (2012), no.~2, 949--976.

\bibitem[BT]{BT} H.~Bachmann and K.~Tasaka, \textit{The double shuffle relations for multiple Eisenstein series}, Nagoya Math. J. \textbf{230} (2018), 180--212.

\bibitem[G]{G} A.~B.~Goncharov, \textit{Galois symmetries of fundamental groupoids and noncommutative geometry}, Duke Math. J. \textbf{128} (2005), no.~2, 209--284.

\bibitem[GKZ]{GKZ} H.~Gangl, M.~Kaneko and D.~Zagier, \textit{Double zeta values and modular forms}, in: Automorphic Forms and Zeta Functions, World Scientific, 2006, 71--106.

\bibitem[Y]{Y} J.~Yu, \textit{Modularity criteria for multiple Eisenstein series}, preprint, \href{https://arxiv.org/abs/2610.21291}{arXiv:2610.21291}, 2026.

\bibitem[Z]{Z} D.~Zagier, \textit{Evaluation of the multiple zeta values $\zeta(2,\ldots,2,3,2,\ldots,2)$}, Ann. of Math. (2) \textbf{175} (2012), no.~2, 977--1000.
\end{thebibliography}
\end{document}